\documentclass[11pt,letterpaper,journal]{IEEEtran}

\usepackage{setspace}
\onecolumn

\newcommand{\superscript}[1]{\ensuremath{^{\textrm{#1}}}}

\usepackage{cite}

\ifCLASSINFOpdf
 \usepackage[pdftex]{graphicx}

\else
\fi
\usepackage[cmex10]{amsmath}
\ifCLASSOPTIONcompsoc
\usepackage[caption=false,font=normalsize,labelfont=sf,textfont=sf]{subfig}
\else
\usepackage[caption=false,font=footnotesize]{subfig}
\fi

\begin{document}
\title{Aggregated Demand Modelling Including Distributed Generation, Storage and Demand Response}

\author{\IEEEauthorblockN{Hesamoddin Marzooghi\IEEEauthorrefmark{1}, David J. Hill\IEEEauthorrefmark{1}\IEEEauthorrefmark{2}, and Gregor Verbi\v{c}\IEEEauthorrefmark{1}\\}
\IEEEauthorblockA{\IEEEauthorrefmark{1}School of Electrical and Information Engineering, The University of Sydney, Sydney, NSW, Australia,\\}
\IEEEauthorblockA{Emails:\{hesamoddin.marzooghi, david.hill, gregor.verbic\}@sydney.edu.au\\}
\IEEEauthorblockA{\IEEEauthorrefmark{2}Department of Electrical and Electronic Engineering, The University of Hong Kong, Hong Kong,\\}
\IEEEauthorblockA{Email:\ dhill@eee.hku.hk}}


%


\IEEEpeerreviewmaketitle

\maketitle
\begin{abstract}
It is anticipated that penetration of renewable energy sources (RESs) in power systems will increase further in the next decades mainly due to environmental issues. In the long term of several decades, which we refer to in terms of the future grid (FG), balancing between supply and demand will become dependent on demand actions including demand response (DR) and energy storage. So far, FG feasibility studies have not considered these new demand-side developments for modelling future demand. In Australia, installed rooftop photovoltaic (PV) generation has been increasing significantly in recent years, and this is increasingly influence the nett load profile for FG feasibility studies. This paper proposes an aggregate nett load or demand model to be used at higher voltage levels considering the effect of rooftop PV, energy storage and DR for FG scenarios, which is inspired by the smart home concept. The proposed model is formulated as an optimization problem aiming at minimizing the electricity cost. As a case study, the effect of the load model is studied on the load profile, balancing and loadability of the Australian National Electricity Market in 2020 with the increased penetration of wind and solar generation.
\end{abstract}


\section{Introduction}

\par Australia's renewables portfolio target aims at increasing the penetration of renewable energy sources (RESs) to 20\% and 50\% by 2020 and 2050, respectively \cite{DepartofResources}. As a result, maintaining the balance between supply and demand will become more dependent on the demand-side to achieve balance with fluctuating renewable power generation. In conventional power systems, the generation is able to follow demand by a combination of dispatch and regulation processes. However, in future grids (FGs) both demand and generation are variable and actions such as demand response (DR) alongside energy storage will be needed to achieve balance. So, the paradigm is changing from generation following demand to demand following generation. This leads to a new modelling question, namely how to represent the aggregated nett demand or load (including distributed generation (DG), storage and DR) for modelling in FG scenarios.

\par In particular, installed capacity of rooftop photovoltaic (PV) has been increasing considerably in the world, and greater penetration of battery storage is anticipated \cite{EPRI, RockyMountain, PV, AEMO2012A, AEMO2012, CSIRO, ATA }. Global installed capacity of PV has increased from approximately 4 GW in 2003 to nearly 128 GW in 2013 \cite{EPRI}. Such a deployment of PV is due to electricity price increases, government incentives and also worldwide drop of PV capital costs \cite{EPRI, RockyMountain}. The Electric Power Research Institute (EPRI) has recently shown that users equipped with PV and storage will be integrated into the USA grids, and provide responsive demand in the future \cite{EPRI}. Another study using a copper plate transmission model has demonstrated that users equipped with PV-plus-battery systems will reach retail price parity from 2020 in the USA grids \cite{RockyMountain}. Also, it can be argued that demand management would accelerate this phenomenon. Until 2013, over 3 GW of PV generation had been installed in Australia mostly by residential and commercial customers \cite{PV}. Moreover, installed battery storage capacity has been increasing in the past couple of years \cite{AEMO2012A, AEMO2012}. It is predicted that PV and storage penetration in Australia will increase further in the next decades \cite{AEMO2012A, AEMO2012, CSIRO}, and, this will greatly influence the load profile and future direction of Australia's electricity system \cite{CSIRO}. The Alternative Technology Association (ATA) in Australia has also suggested that due to high penetration of PV-plus-storage system in the Australian National Electricity Market (NEM), retail price parity appears highly plausible by 2020 \cite{ATA}. The Australian Energy Market Operator (AEMO) has also noted an increased difficulty of predicting the future demand because of the rapid growth of those sources in the NEM \cite{AEMO2012A}. Furthermore, it is clear that for studying FG scenarios of uptake of the various demand technologies (DGs, storage, DR), models which can represent different levels of each technology are needed.

\par The literature review below considers studies which relate directly to meeting new modelling and analysis challenges in FGs \cite{Energy2010, Elliston2012, Elliston2013, Budischak2013, Hart2011}. A preliminary study by the University of Melbourne Energy Research Institute has proposed an electrical grid for Australia in 2020 relying 100\% on RESs \cite{Energy2010}. However, balancing and performance of the proposed grid was not studied which makes the proposal highly speculative. The University of New South Wales researchers have analysed the feasibility of 100\% RES scenarios considering a copper plate transmission model for the NEM \cite{Elliston2012, Elliston2013}. They have demonstrated that the NEM can be powered most of the time entirely on RESs within the specified NEM reliability standard. Also, they have determined the least-cost mix of 100\% RESs scenario including wind farms (WFs), utility PV, concentrated solar plants (CSPs) with thermal storage, hydro and biofuelled gas turbines (GTs). Beside national FG feasibility studies, some studies are reported for other electrical grids in the world. USA researchers have proposed a combination of RESs (i.e. onshore and offshore WFs, utility PV and fuel cells) and conventional generation (i.e. GTs) for the future of the PJM network using a copper plate network \cite{Budischak2013}. They have illustrated that demand in the PJM network can be met 90-99.9\% of the time with the suggested mix of RESs, while RESs are predicted to be at price parity in 2030 with the existing PJM generation system. In a similar study by other USA researchers, the least-cost mix of RESs (i.e. WFs, CSPs, utility PV, hydro and geothermal) and conventional generation (i.e. GTs) has determined for California in 2050\cite{Hart2011}.

\par The existing studies have demonstrated that relying on higher penetration of RESs is possible, and fossil fuel use can be limited in the future. However, all of the above studies have used conventional load models and so neglected the influence of newer demand-side technologies for modelling nett future demand. Due to the significant effect of loads on performance and stability of power systems, it can be expected that this will affect the results of FG feasibility studies significantly. The Australian Commonwealth Scientific and Industrial Research Organisation (CSIRO) FG Forum has debated that demand management will play a major role in the FG of Australia \cite{CSIRO}. They have also proposed a simple aggregated demand model at state levels in the NEM. However, a generic modelling approach is still required which can be used for any granularity level in the grid (e.g. from a city, a state or even the whole network). Furthermore, all of the reviewed studies have focused on balancing by using a simplified grid models such as the copper plate model. These assumptions can influence the results of power system studies significantly.
 
\par In this study, we make a step towards aggregate demand modelling suitable for FG scenarios. The proposed load model considers aggregated demand-side PV and storage, and is inspired by the smart home concept \cite{Tischer2011}. The load model is formulated as an optimization problem aiming at minimizing the customer electricity costs. The electricity price is predicted by core vector regression (CVR) \cite{Tsang2006} using a combination of historical market data from the AEMO and simulated market prices for the future. Then, as a case study, the effect of the model on the load profile, balancing and loadability of the NEM in 2020 is studied using a 14-generator model \cite{Gibbard2010}. The electricity market model is built in PLEXOS based on the suitably modified 14-generator model. Five scenarios are analysed with one business as usual (BAU), and four different levels of DR with renewable integration. (i) For the BAU Scenario in 2020, the electricity supply is dominated by coal, gas, hydro and biomass; and (ii) for the Renewable Scenarios, some of the conventional coal generators in Queensland and South Australia are replaced with CSPs with storage and WFs respectively, as suggested in \cite{DepartofResources, AEMO2012A, Energy2010} to meet Australia's RES target. The Renewable Scenarios are analysed with both the conventional and the proposed load models. Simulation results show that demand management can improve balancing and decrease the required energy from the backup supply with the increased intermittent supply in the grid.

\par The remainder of the paper is organised as follows: Section II proposes the aggregated load model considering demand-side technologies. Section III describes the test-bed and the electric market assumptions and modelling. Section IV describes simulation scenarios, and discusses simulation results. Section V consists of the conclusion of the simulations.

\section{Aggregated Load Model Considering Demand-Side Technologies}

\par Aggregated load models are commonly used in system studies to reflect the combined effect of numerous physical loads \cite{Kundur, Van}. These can be inspired by a physical devices, e.g. using a large induction motor to represent all the motors connected, or by data-driven approaches. This paper proposes an aggregate load model for studies at higher voltage levels inspired by the smart home concept \cite{Tischer2011}. A smart home can be thought of as an automated residential building that uses its own DGs for robustly managing energy consumption to reduce electricity costs. A smart home energy management system (SHEMS) implements an algorithm to schedule DGs and storage, and so achieve demand-side control \cite{Tischer2011}. The proposed model includes aggregated PV, storage and DR assuming users are price takers, i.e. the effect of their actions is not reflected in the electricity price. This assumption is usually considered when the number of users is large, and the amount of information provided to each user is limited \cite{Samadi2012}.

\subsection{Optimization model}

\par The proposed load model aims at minimizing the electricity cost for end-users by reducing power consumption from the electrical grid. This aim can be achieved by controlling the battery state of charge (SOC) at the beginning of those hours. Each decision horizon for the model (i.e. 24-hour period) is divided into one hour time-steps, giving a total of $H$=24 time-steps; denote a particular time-step by $h$. The objective function of the model, which is solved using the linear programming approach, can be written as:

\begin{subequations}
     \begin{align}
         & \underset{}{\text{min}}
         && \sum_{h=1}^{H} C_{el}(h)P_{g}(h),\\
         & \text{s.t.}
         & & P^{min}_{grid}\leq P_{g}(h)\leq P^{max}_{grid},\\ 
         &
         & & B^{dis}_{rate}\leq P_{b}(h)\leq B^{cha}_{rate},\\ 
         &
         & & B_{SOC}(1) = B^{min}_{SOC},\\
         &
         & & B_{SOC}(h) = B_{SOC}(h-1) + P_{b}(h-1)\quad\forall h>1,\\
         &
         & & B^{min}_{SOC}\leq B_{SOC}(h)\leq B^{max}_{SOC},\\
         &
         & & P_{g}(h) = P_{L}(h)+ \eta P_{b}(h)-P_{pv}(h),
         \end{align}
\end{subequations}
where, $C_{el}$ and $P_{g}$ are the electricity price and the grid power, respectively. In the (1b), the grid power is constrained by the maximum, $P^{max}_{grid}$, and the minimum, $P^{min}_{grid}$, electrical grid power. If, $P_{g}(h)\geq0$ the power is flowing from the grid to the demand. Otherwise, the power is flowing back to the grid. Constraints (1c) to (1f) represent battery storage power, $P_{b}$, and its SOC limit, where, $B^{cha}_{rate}$, $B^{dis}_{rate}$, $B_{SOC}$, $B^{max}_{SOC}$ and $B^{min}_{SOC}$ are charging and discharging rate, SOC, maximum and minimum SOC, respectively. If $P_{b}(h)\geq0$, battery stores energy, and if $P_{b}(h)\leq0$, it discharges the stored energy. The last constraint (1g) is the balance equation for the load model, where, $P_{pv}$, $P_{L}$ and $\eta$, are the aggregated PV and demand power, and the battery efficiency, respectively. 

\par The model in (1a)-(1g) augments a conventional balancing load model by the aggregated effect of numerous price taker users equipped with PV-plus-storage systems scheduled with the SHEMS. As the model is formulated, price responsive users shift their consumption from expensive time-slots to cheaper ones to utilize cheaper electricity produced by DGs. Also, one of the advantages of the load model is that it has the capability to be used for any granularity level in the grid (e.g. from a busbar to a city, a state or even the whole grid).
 
\par Considering price taker entities is a realistic assumption for a large number of users \cite{Samadi2012}. However, high penetration of those users might affect power system performance, as discussed in more detail in Section IV. To solve the above optimization problem, the electricity price, aggregated demand and PV power are required for each time-step. The next subsections describe how these variables are obtained for the model.

\subsection{Electricity price prediction}

\par In the model (1a)-(1g), the electricity price, $C_{el}$, plays a key role. Amongst different techniques, we chose intelligent methods that have been used successfully in the past for electricity price prediction \cite{Szkuta1999, Yamin2004, Niu2010}. In this study, we used CVR \cite{Tsang2006} for the electricity price prediction which requires a short training time and small memory requirement for a large training dataset. The modelling algorithm and inputs for the CVR were inspired by the models reported in \cite{Szkuta1999, Yamin2004}. So, the interstate line limits, predicted demand, time factor, capacity, type and area of generators are chosen as the inputs and the electricity price is selected as the output for the CVR. The electricity price predictor has been trained using a combination of historical data from the AEMO and simulated market prices for the future (For the future, the electricity prices have been simulated in PLEXOS to cater for situations where part of the generation bids in to the market at zero cost). Using the trained CVR, the electricity price is predicted for 2020, and is given to users.

\subsection{PV and battery storage}

\par The aggregated demand ($P_{L}$) and PV ($P_{pv}$) power are uncertain variables in a particular time-step of the load model. The AEMO has proposed 16 zones for the NEM to capture differences in generation technology capabilities, costs, weather and so on in the future \cite{AEMO2012A}. In this study, hourly demand and PV power are obtained from the AEMO predications for 2020. The PV output power is reported for 16 zones, however, the reported demand data is aggregated across each region of the NEM (i.e. Queensland (QLD), New South Wales (NSW), Victoria (VIC) and South Australia (SA)). The predicted demand data is divided between the zones based on the available loads of the 14-generator model in each zone and their default values. The 14-generator model of the NEM and the match between 14-generator model and the AEMO zones are described in Section III. 

\par According to the AEMO, the percentage of the residential and commercial customers with PV will be 23\%, 28\% and 36\% for low, moderate and high uptake scenarios in 2020, respectively \cite{AEMO2012}. In this paper, price-responsive users are chosen from the residential and commercial customers, which are considered to account for 60\% of the total system load in the NEM in 2020 \cite{AEMO2012A}. The industrial customers are left unaffected. Also, the percentage of the residential and commercial customers with PV are considered 20\%, 30\% and 40\% for low, medium and high uptake scenarios, respectively. Table~\ref{tab:The aggregated storage and PV capacities for each region of the NEM for different DR scenarios} shows the aggregated storage and PV capacities for each region of the NEM for different uptake scenarios. The chosen PV capacities for each region of the NEM for different uptake scenarios are inspired by the AEMO study \cite{AEMO2012}.  Also, the chosen battery storage capacities roughly correspond with a typical PV and storage capacity for a household in Australia (i.e. 3kW PV and 10kWh battery storage).

\begin{table}  
\centering
\caption{the aggregated storage and pv capacities for each region of the nem for different uptake scenarios}
\begin{tabular}{| c | c | c | c |}\hline
Region & Scenario & $B^{min}_{SOC}$-$B^{max}_{SOC}$ & PV capacity \\
       &          &         (GWh)                   &   (GW)      \\\hline
    & Low    & 0.4-4.3  & 1.3 \\
QLD & Medium & 0.6-6.4  & 1.9 \\
    & High   & 0.9-8.5  & 2.6 \\\hline
    & Low    & 0.7-6.7  & 2.0 \\
NSW & Medium & 1.0-10.1 & 3.0 \\
    & High   & 1.4-13.5 & 4.1 \\\hline
    & Low    & 0.5-5.0  & 1.5 \\
VIC & Medium & 0.8-7.5  & 2.3 \\
    & High   & 1.0-10.0 & 3.0 \\\hline
    & Low    & 0.1-1.2  & 0.3 \\
SA  & Medium & 0.2-1.7  & 0.5 \\
    & High   & 0.2-2.3  & 0.7 \\\hline
    & Low    & 1.7-17.0 & 5.0 \\
NEM & Medium & 2.5-25.0 & 7.5 \\
    & High   & 3.4-34.0 & 10.5 \\\hline        
\end{tabular}
\label{tab:The aggregated storage and PV capacities for each region of the NEM for different DR scenarios}
\end{table}

\section{The Australian NEM Model}

\par The main contribution of this paper is to propose the load or demand model considering new technologies (DG, storage and DR). In Section IV, the effect of the model will be illustrated on the load profile and performance of the NEM in 2020. A 14-generator model of the NEM, which was originally proposed for small-signal stability studies \cite{Gibbard2010} is used as the test-bed, where market effects and constraints are also taken into account.

\subsection{Test-bed assumptions and modelling}

\par
\begin{figure}  
\centering
\includegraphics[width=11cm, height=16cm]{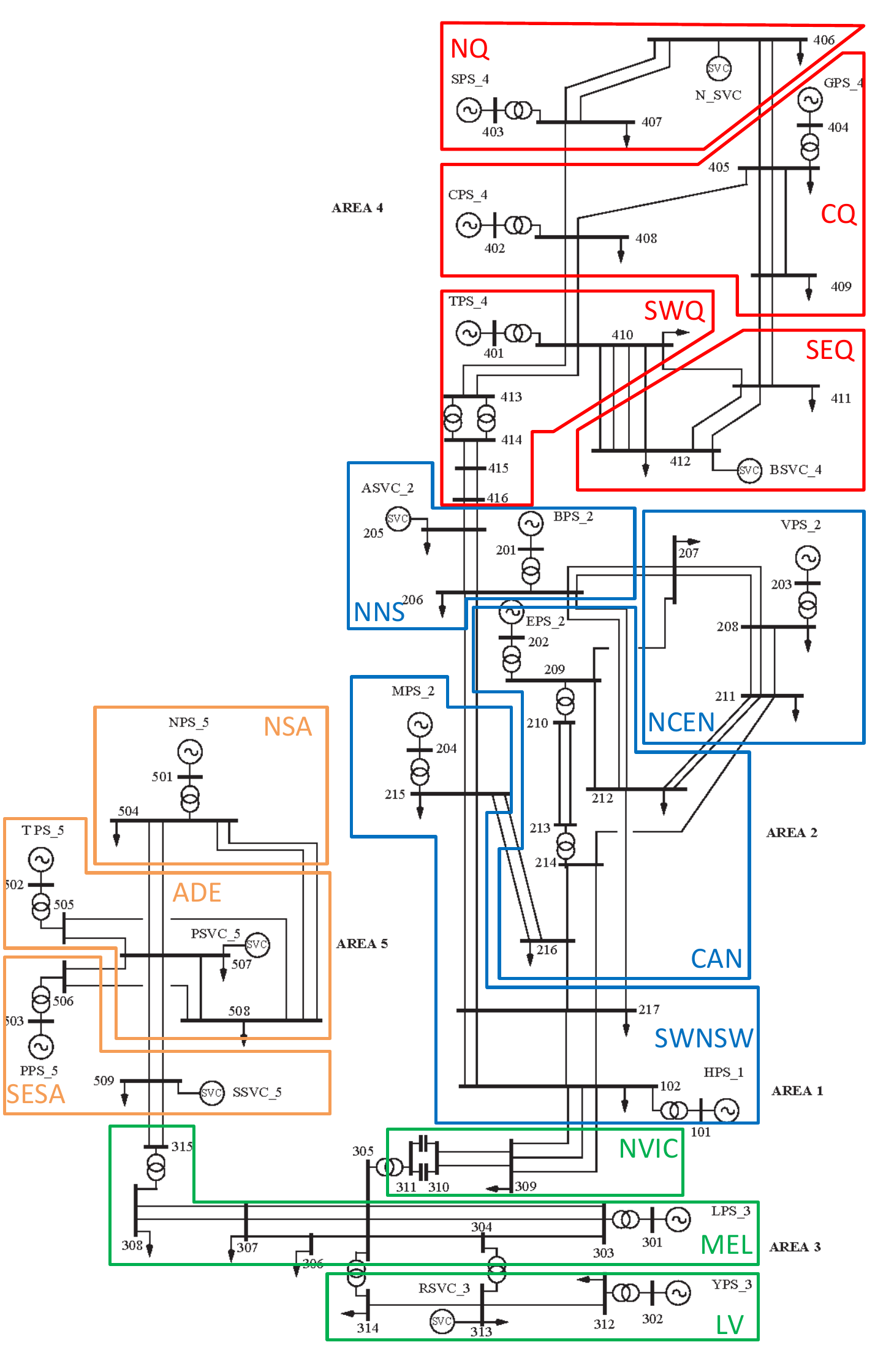}
\caption{14-generator model of the NEM}
\label{figure:14-generator model of the NEM}
\end{figure}

 The schematic diagram of the 14-generator model of the NEM is shown in Figure~\ref{figure:14-generator model of the NEM}. Areas 1 to 5 represent Snowy Hydro (SH), NSW, VIC, QLD and SA, respectively. In order to extract data for the load model and generators in 2020, the 14-generator model of the NEM is matched with the 16 zones according to the AEMO's planning document \cite{AEMO2012A}, as shown in Figure~\ref{figure:14-generator model of the NEM}. The modified 14-generator model of the NEM is then modelled in PLEXOS and MATLAB (MATPOWER) for the market simulations, balancing and loadability studies, respectively.

\subsection{Electricity market assumptions and modelling}

\par The market simulations of this paper are done in PLEXOS, following the dispatch process used by the AEMO. The resolution of the market simulations is taken as one hour. Combinations of coal, GT, hydro and biomass are considered for the NEM to supply the load in 2020 in the BAU Scenario. The generator technologies in this study are assumed as: pulverized supercritical coal generators based on black and brown coal, combined cycle GT based on natural gas, solar parabolic trough CSP with solar multiple of 2 and 12 hours thermal storage, biomass based on landfill gas (LFG) and WFs. The coal power plant in NSW and QLD are assumed to be black coal-based, and in VIC and SA are assumed to be brown coal-based \cite{DepartofResources, AEMO2012A}. For each hour of the year, the mixed integer linear solver dispatches the generation, in merit order, to meet demand for that hour. In market simulations, the fossil-fuel generators were assumed to bid at their respective short-run marginal costs (SRMC), calculated using the predicted fuel price, thermal efficiency and variable O\&M \cite{AEMO2012A, Energy2012}, while the SRMC of renewable generation is assumed to be zero. Table~\ref{tab:The predicted cost of the chosen technologies in 2020} lists the SRMC of the generators in 2020. The maximum and the minimum interstate line flow between NSW/QLD, NSW/VIC and VIC/SA are considered to be (600, -1000), (500, -1500) and (500, -500) MW, respectively. The interstate line flows roughly correspond with the NEM limits \cite{AEMO2012A}. The market model also considers the minimum stable level of generators reported in \cite{AEMO2012A}.

\begin{table}
\centering
\caption{the srmc of the generators in 2020}
\begin{tabular}{| c | c |  c | c | c |}\hline
Generator  & Type           & AEMO zone  &  SRMC\\
           &                 &           & (\$/MWh)\\\hline
 BPS\_2    &  Black coal     & NNS       & 28.45  \\
 EPS\_2    &  GT             & CAN       & 69.20  \\
 MPS\_2    &  Black coal     & SWNSW     & 27.43  \\
 VPS\_2    &  Black coal     & NCEN      & 26.40  \\
 LPS\_3    &  Biomass        & MEL       & 39.50  \\
 YPS\_3    &  Brown coal     & LV        & 21.88  \\
 CPS\_4    &  Black coal     & CQ        & 26.14  \\
 GPS\_4    &  Black coal     & CQ        & 26.14  \\
 SPS\_4    &  Black coal     & NQ        & 32.74  \\
 TPS\_4    &  GT             & SWQ       & 73.84  \\
 NPS\_5    &  Brown coal     & NSA       & 30.89  \\
 PPS\_5    &  Brown coal     & SESA      & 30.89  \\
 TPS\_5    &  GT             & ADE       & 69.20  \\\hline         
\end{tabular}
\label{tab:The predicted cost of the chosen technologies in 2020}
\end{table}

\par
\begin{figure}
\centering
\includegraphics [width=12.5cm, height=10.5cm]{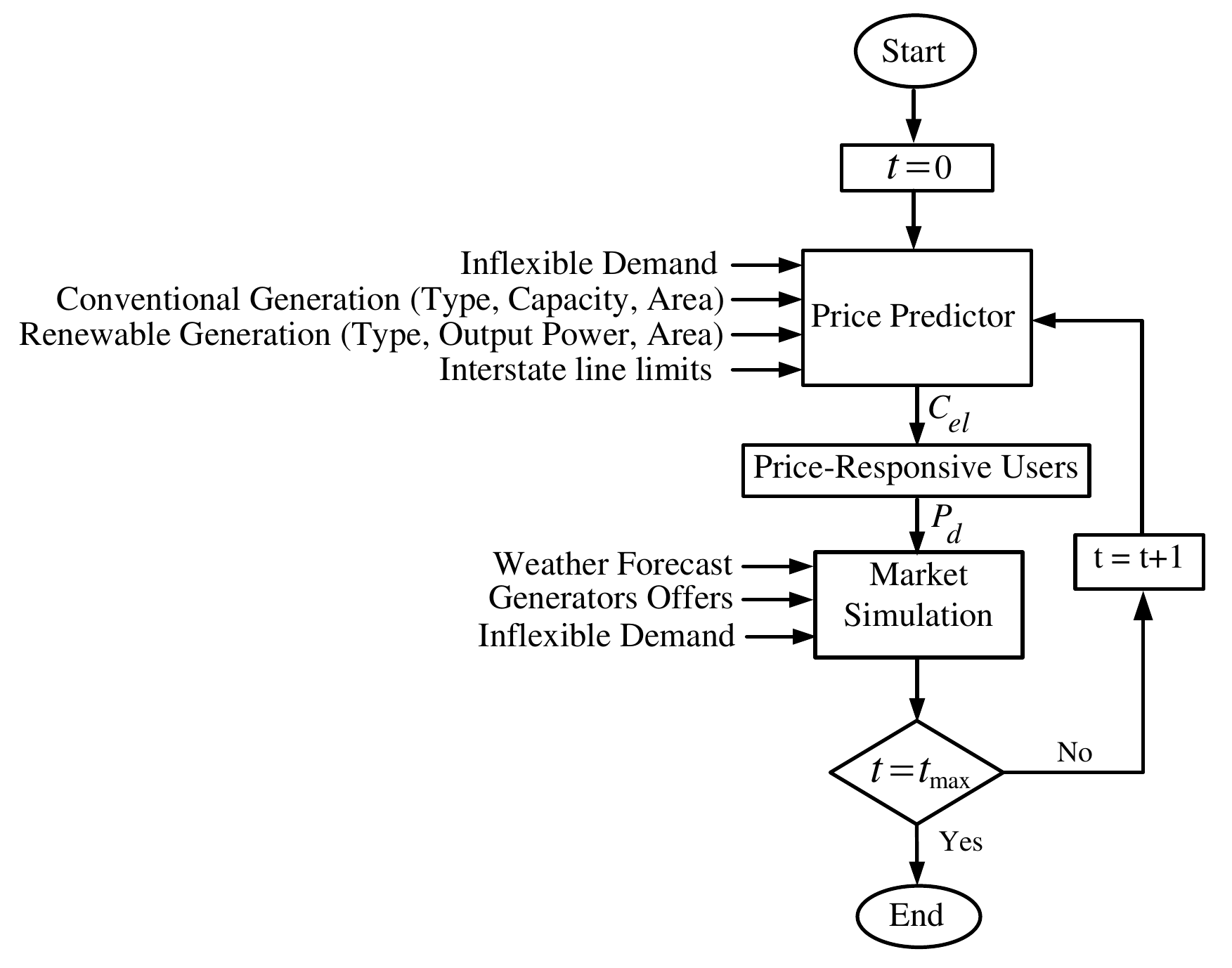}  
\caption{The relation between the electricity price, demand and the electricity market}
\label{figure:DR modeling flowchart}
\end{figure}

 Figure~\ref{figure:DR modeling flowchart} illustrates the relation between the electricity price, demand and the electricity market. The CVR-based model predicts the electricity price signal for each hour, and price-responsive users respond to the signal. The demand data is then communicated to the electricity market to be used for the dispatch process. If supply cannot meet the demand, the hour is recorded as the unserved hour. However, if available generation exceeds demand (i.e. due to high generation of intermittent RESs), the surplus power is recorded as dumped energy and that hour is marked as a dumped hour. Finally, the dispatch results are used as input for the MATPOWER for balancing and loadability studies.

\section{Simulation Scenarios and Results}

\par The following subsections describe simulation scenarios, and demonstrate the effect of the proposed model on the load profile, balancing and loadability of the NEM in 2020 with increased penetration of WFs and CSPs. 

\subsection{Simulation scenarios}

\par Five scenarios will be analysed with one BAU and four different levels of DR with renewable integration. For the BAU Scenario, a combination of coal, gas, hydro and biomass are considered to supply the load in 2020 (i.e. Scenario 1). Then, some of the conventional coal generators in QLD and SA are replaced with CSPs together with storage and WFs respectively to meet Australia's RES target. Displacement of the conventional generators in the Renewable Scenarios and the chosen capacity for the RESs are inspired by the studies in \cite{AEMO2012A, Energy2010}. NPS\_5 in SA is replaced with a WF with the capacity of 3 GW using NSA data. Also, SPS\_4 and GPS\_4 in QLD are replaced with two CSPs with the capacity of 4.5 GW each and using NQ and CQ data, respectively. The operating strategy for CSPs has been determined using PLEXOS. It was found that delaying CSP output by 12 hours minimizes the unserved and dumped energy. The RESs serve about 20\% of the total demand energy in the Renewable Scenarios. The Renewable Scenarios are evaluated with the conventional and the proposed load to study the effect of demand-side technologies on the system performance with the increased penetration of RESs. In the rest of the paper, the Renewable Scenarios with conventional load, low, medium and high uptake of demand-side technologies are called Scenarios 2 to 5, respectively.

\subsection{Load profile}

\par
\begin{figure} 
{\subfloat[]{\includegraphics[width=9.4cm, height=6.5cm]{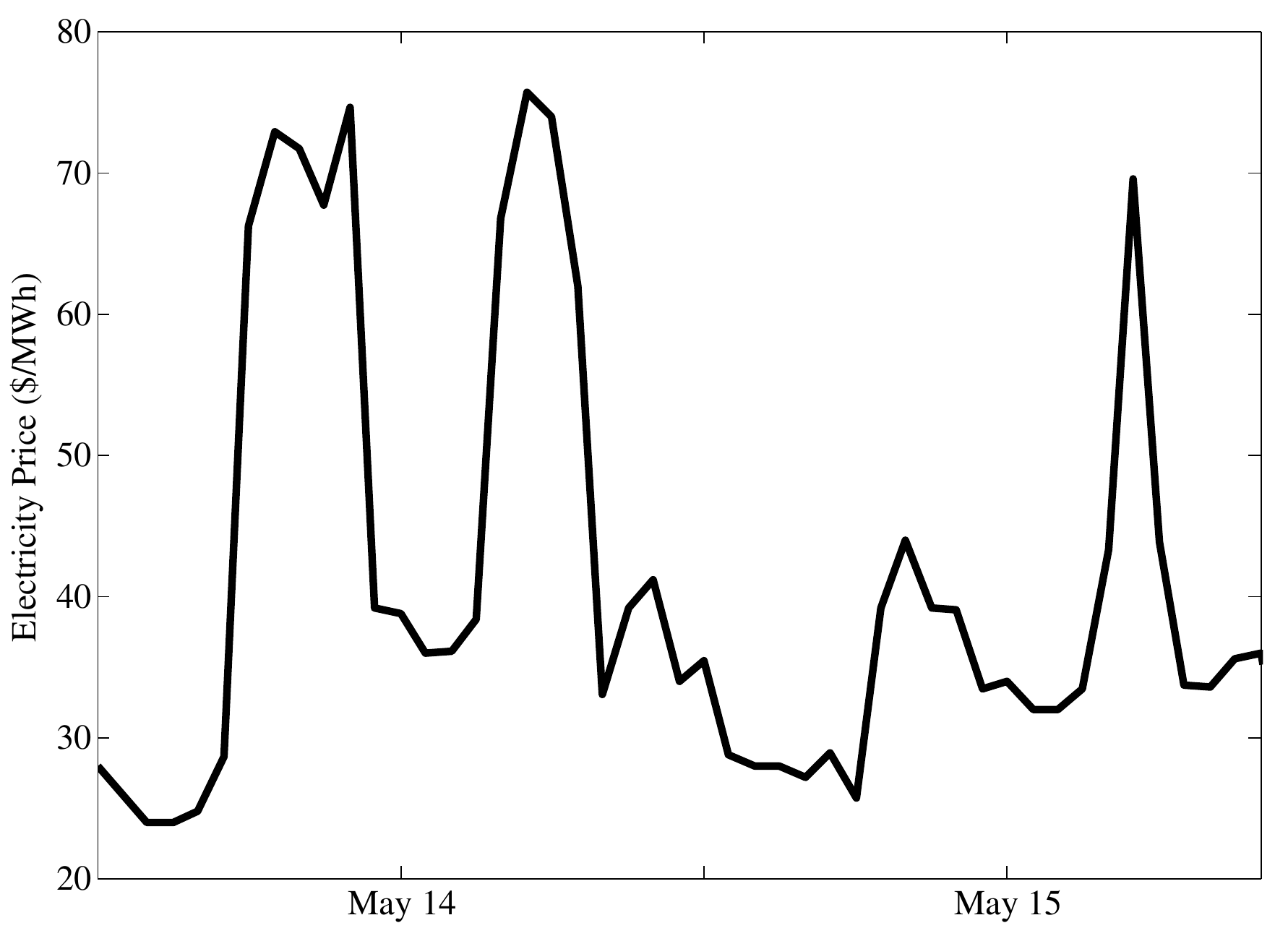}
\label{fig:first_caseA1}}}
{\subfloat[]{\includegraphics[width=9.4cm, height=6.5cm]{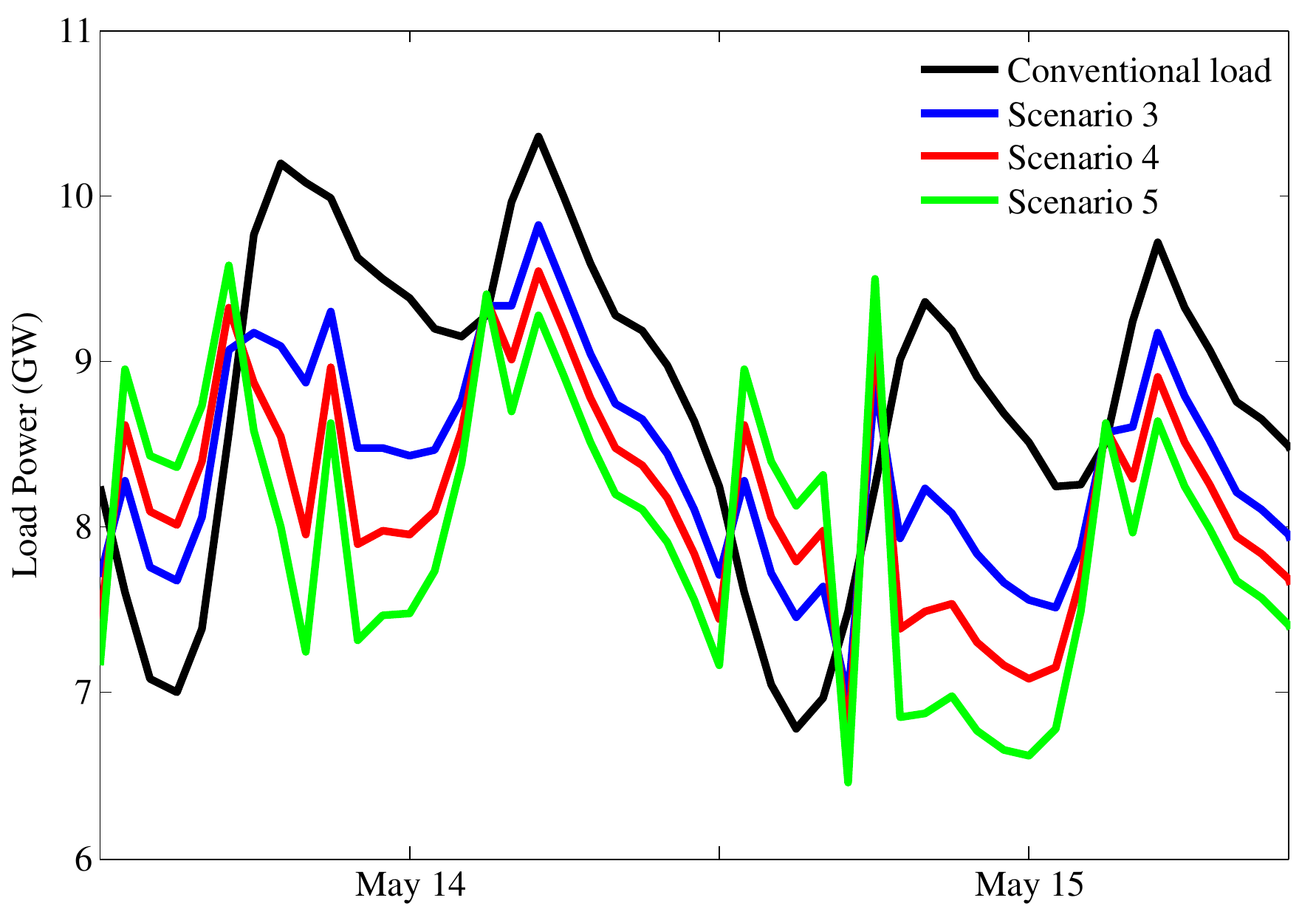}
\label{fig_second_caseA2}}} 
{\subfloat[]{\includegraphics[width=9.4cm, height=6.5cm]{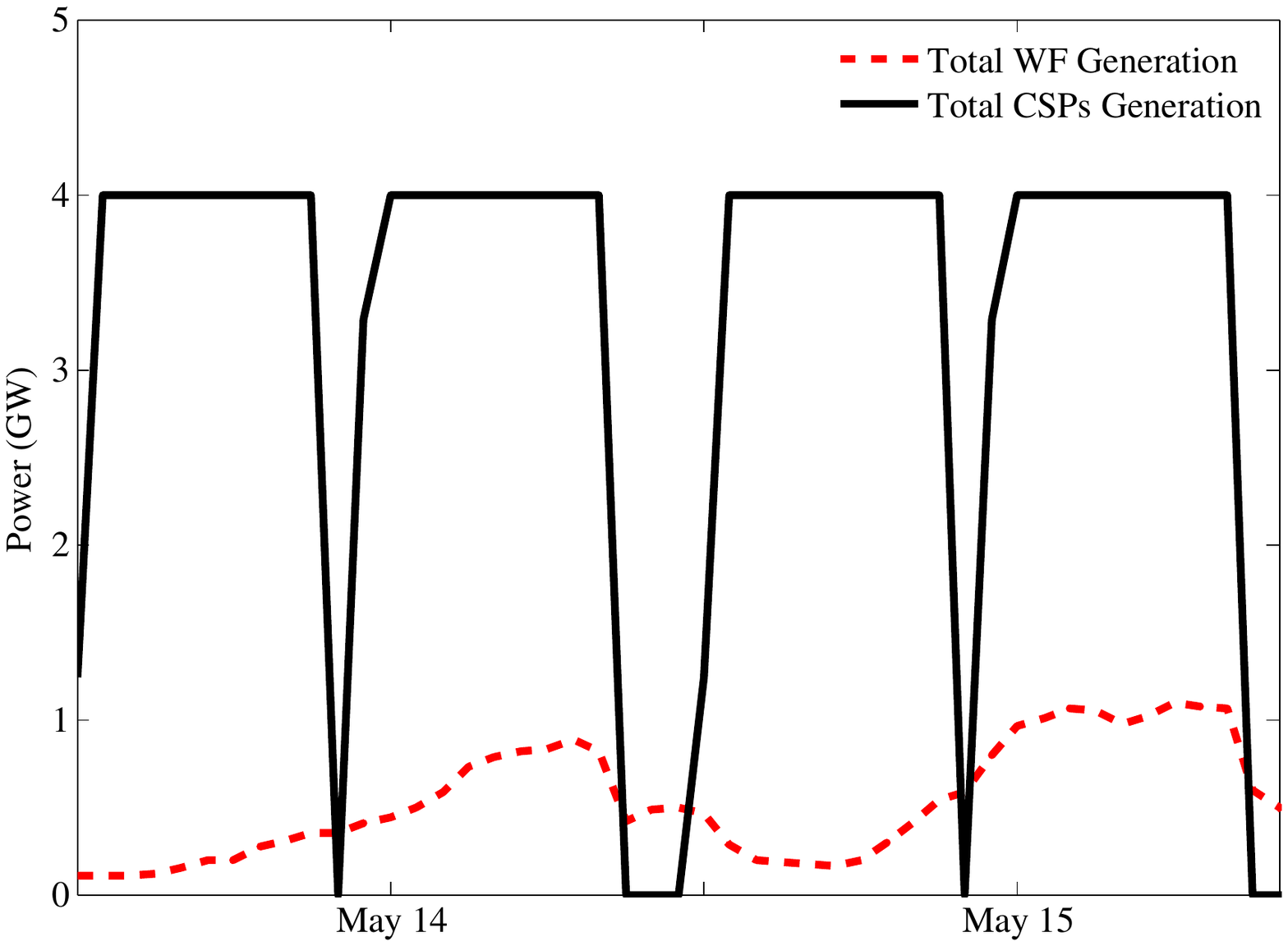}
\label{fig_second_caseA2}}} 
{\subfloat[]{\includegraphics[width=9.4cm, height=6.5cm]{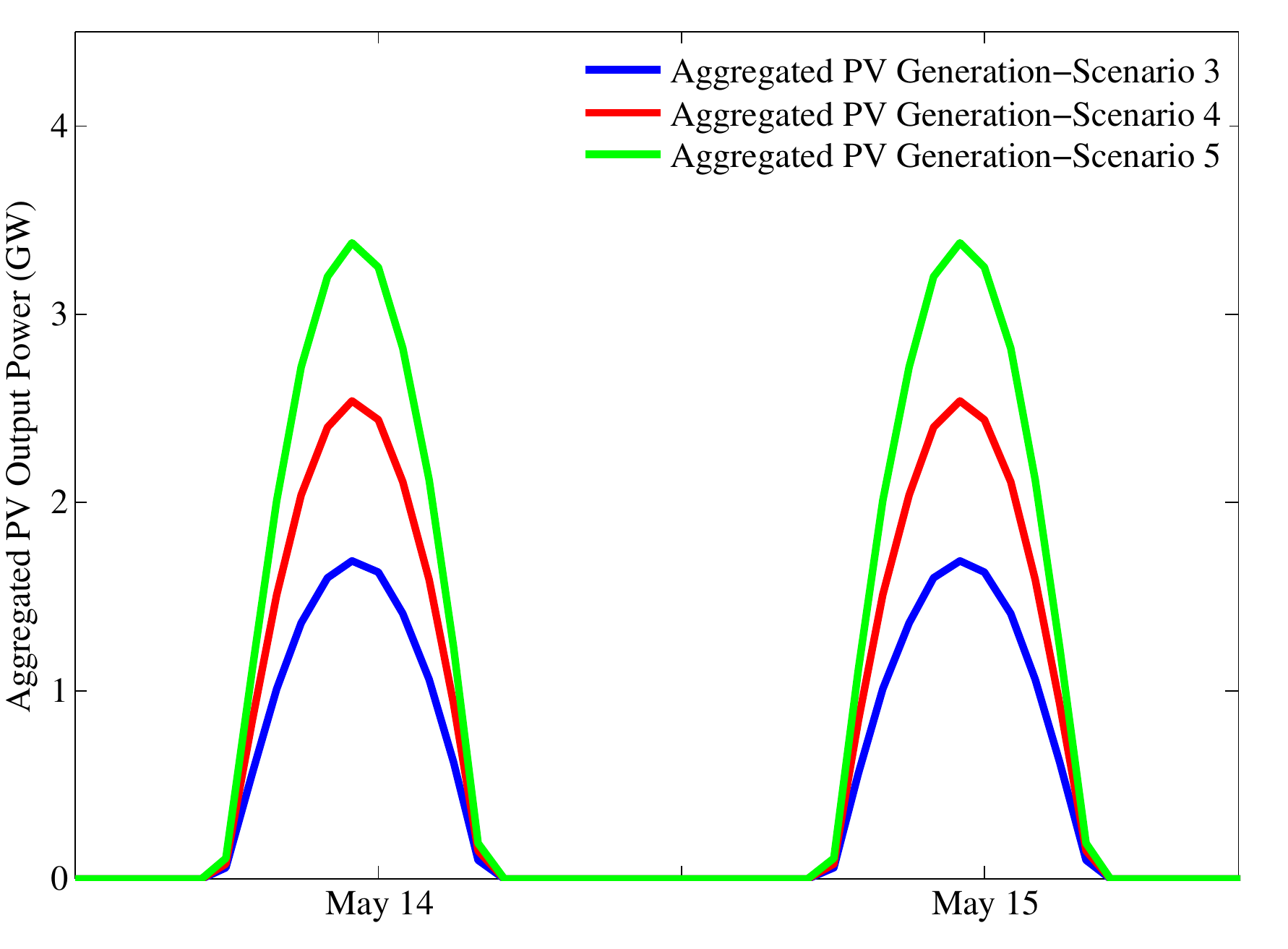}
\label{fig_second_caseA2}}} 
\caption{(a) Electricity price, (b) demand profile, (c) total RESs generation, and (d) aggregated PV generation for Scenarios 2 to 5 in NSW during the 14\superscript{th} and the 15\superscript{th} of May 2020}
\label{figure:NSW}
\end{figure}

 The load in Section II is modelled using the predicted electricity price and the aggregated demand and PV power. Figure~\ref{figure:NSW} (b) shows the effect of different PV-plus-storage uptake scenarios on the load profile resulting from solving the proposed model in NSW during the 14\superscript{th} and 15\superscript{th} of May 2020. The electricity price signal (Figure~\ref{figure:NSW} (a)) is predicted by the CVR using demand data, RESs generation (Figure~\ref{figure:NSW} (c)), conventional generation availability in the grid, and line flow limits. As it can be seen in Figure~\ref{figure:NSW} (a), when generation from RESs decreases, electricity price mainly increases because conventional generators need to compensate for the lack of generation in the grid. Price-responsive users respond to the electricity price signal and shift their consumption from expensive time-slots to the cheaper ones (i.e. using their PV-plus-storage system) to utilize zero cost electricity produced by RESs more, as it can be seen in Figure~\ref{figure:NSW} (b). This clearly shows that DR can help balance intermittent RESs power and demand in FGs. However, as users are price takers, they may shift their consumption to cheaper time-slots all together which can result in secondary peaks for the system under high penetration of demand-side sources and storage, as shown in Figure~\ref{figure:NSW} (b). Furthermore, management of price-responsive users using this electricity price signal might result in less smooth load profile in comparison with the conventional load profile because the effect of users' action is not reflected in the electricity price signal. To address this, demand response aggregators will likely emerge in the future, which will require the loads to be treated as price anticipator entities, necessitating game-theoretic approaches for the electricity price signal designing \cite{ Samadi2012, Chapman2013}.

\par In the following subsection, both conventional load and the proposed load are used to study the effect of demand-side resources and storage on the balancing and loadability of the NEM with the increased penetration of RESs in 2020.

\subsection{Balancing and loadability results}

\subsubsection{BAU Scenario}
\par
\begin{figure} 
\centering
\includegraphics[width=9.4cm, height=6.5cm]{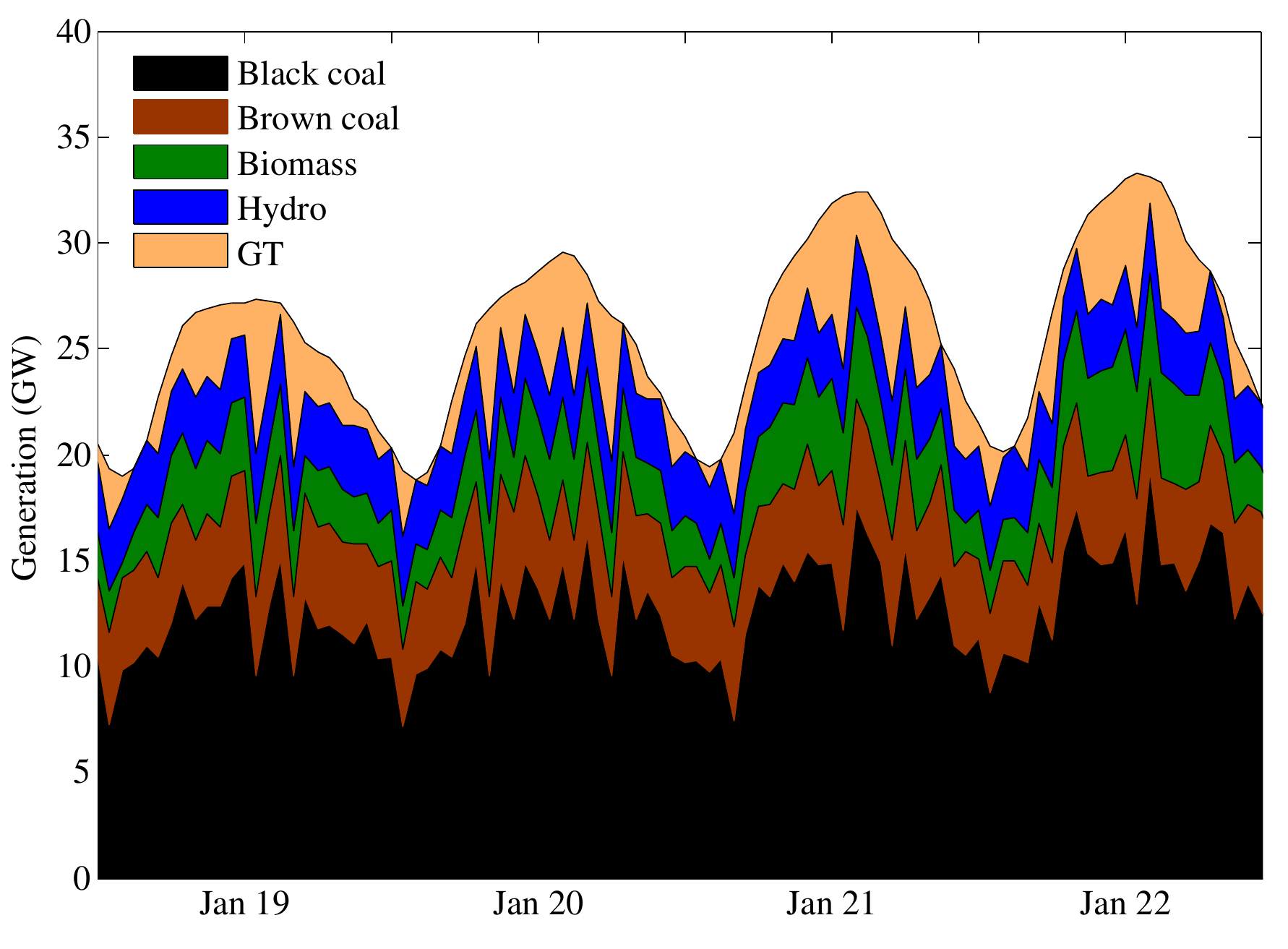}
\caption{Balancing results for the BAU scenario from the 19\superscript{th} to the 22\superscript{nd} of January 2020}
\label{figure:Balancing results for the BAU scenario from 10th to 15th of January 2020}
\end{figure}

 Figure~\ref{figure:Balancing results for the BAU scenario from 10th to 15th of January 2020} shows balancing results for the BAU Scenario in one of the critical summer peaks from the 19\superscript{th} to the 22\superscript{nd} of January 2020. A big portion of the demand is supplied by coal-fired power plants and the peak loads are met with GTs. The results for supplied electrical energy from GTs and average loadability are summarized in Table~\ref{tab:Balancing, loadability and minimum singular value  of $G_{S}$ for scenarios 1 to 4}. It is noteworthy to mention that for loadability calculation, all loads in QLD are assumed to increase uniformly in small steps with constant power factor until power flow fails to converge. Also, it is assumed that all the generators in QLD are scheduled with the same participation factor to pick up the system loads. The loadability is computed for each hour until a step before power flow divergence.

\subsubsection{Renewable Scenarios}

\par Four different PV-plus-storage and DR uptake scenarios were considered for 2020: zero, low, medium and high, denoted Scenarios 2, 3, 4, and 5, respectively in Table~\ref{tab:Balancing, loadability and minimum singular value  of $G_{S}$ for scenarios 1 to 4}. Unserved hours for all scenarios are zero. Comparing the BAU Scenario and the Renewable Scenario with the conventional load (Scenario 2), it can be seen that with the increased penetration of RESs, the loadability is reduced from 27.13 GW to 22.17 GW. Also, the required electrical energy from the backup generation (i.e. GTs) is increased from 18.73 to 18.77 TWh.

\par Compared to Scenario 2, a larger penetration of PV-plus-storage (Scenarios 3 to 5) improves the balancing and loadability, and reduces the required energy from backup supply. The high uptake scenario has the lowest spilled energy and hours followed by medium and low uptake scenarios. The medium uptake scenario (Scenario 4) has, however, the highest loadability. The loadability is increased from 22.17 GW for the Scenario 2  to 25.53 GW for the Scenario 4, which implies a considerable improvement in loadability. Surprisingly, increasing the penetration of PV-plus-storage system beyond a certain point (from medium to high), fails to improve loadability further. This is due to the price-taking assumption of users, so for the high uptake scenario (Scenario 5), secondary load peaks are created as all the users shift consumption to the cheaper time slots. This deteriorates loadability compared to lower uptake scenarios. This clearly shows that with large penetration of demand-side technologies (here PV, storage and DR) and the price taking assumption for the loads, the marginal benefit might become negative. To capture the maximum capability of demand-side resources and storage, demand response aggregators will likely emerge in the future. Demand response aggregators will use game-theoretic approaches for the electricity price signal designing, and, will consider the effect of users' action on the electricity price (i.e. they treat users as price anticipator entities) \cite{ Samadi2012, Chapman2013}. This is part of our ongoing research to reflect the aggregated effect of price anticipator users equipped with PV-plus-storage systems at high voltage levels for FG scenarios.

\begin{table}
\centering
\caption{balancing and loadability results for Scenarios 1 to 5}
\begin{tabular}{| c | c | c | c | c |}\hline
Scenarios & Spilled energy & Spilled hrs & GT energy & Loadability\\
          &  (TWh)         & (\%)        &(TWh)      & (GW)\\\hline
     1    &  -   &   -   & 18.73 & 27.13\\   
     2    & 0.71 & 13.65 & 18.77 & 22.17\\
     3    & 0.66 & 13.03 & 18.33 & 23.78\\
     4    & 0.61 & 12.67 & 17.85 & 25.53\\
     5    & 0.54 & 11.71 & 17.28 & 24.41\\\hline
\end{tabular}
\label{tab:Balancing, loadability and minimum singular value  of $G_{S}$ for scenarios 1 to 4}
\end{table}

\par
\begin{figure} 
{\subfloat[]{\includegraphics[width=9.4cm, height=6.5cm]{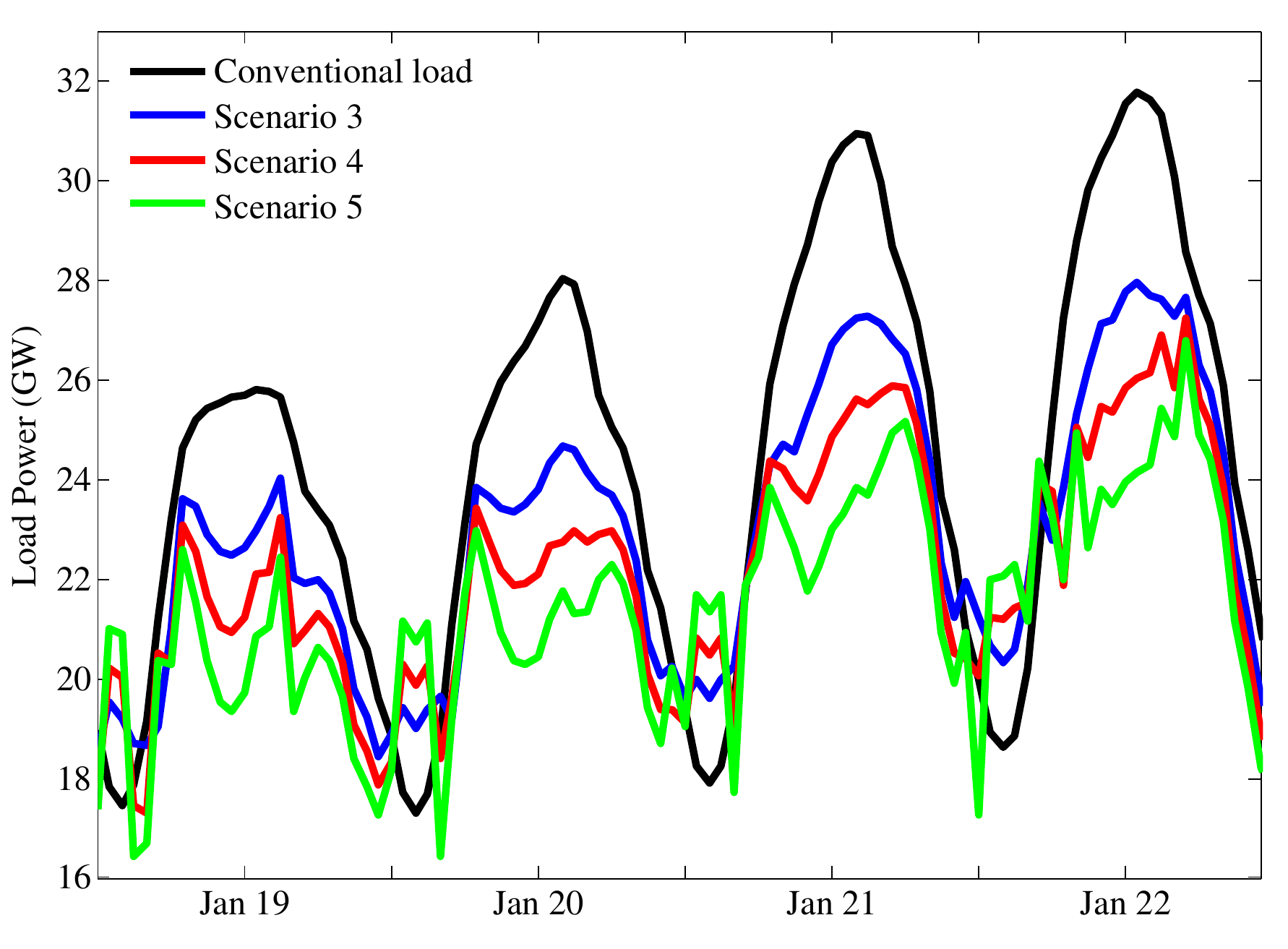}
\label{fig:first_caseA1}}}
{\subfloat[]{\includegraphics[width=9.4cm, height=6.5cm]{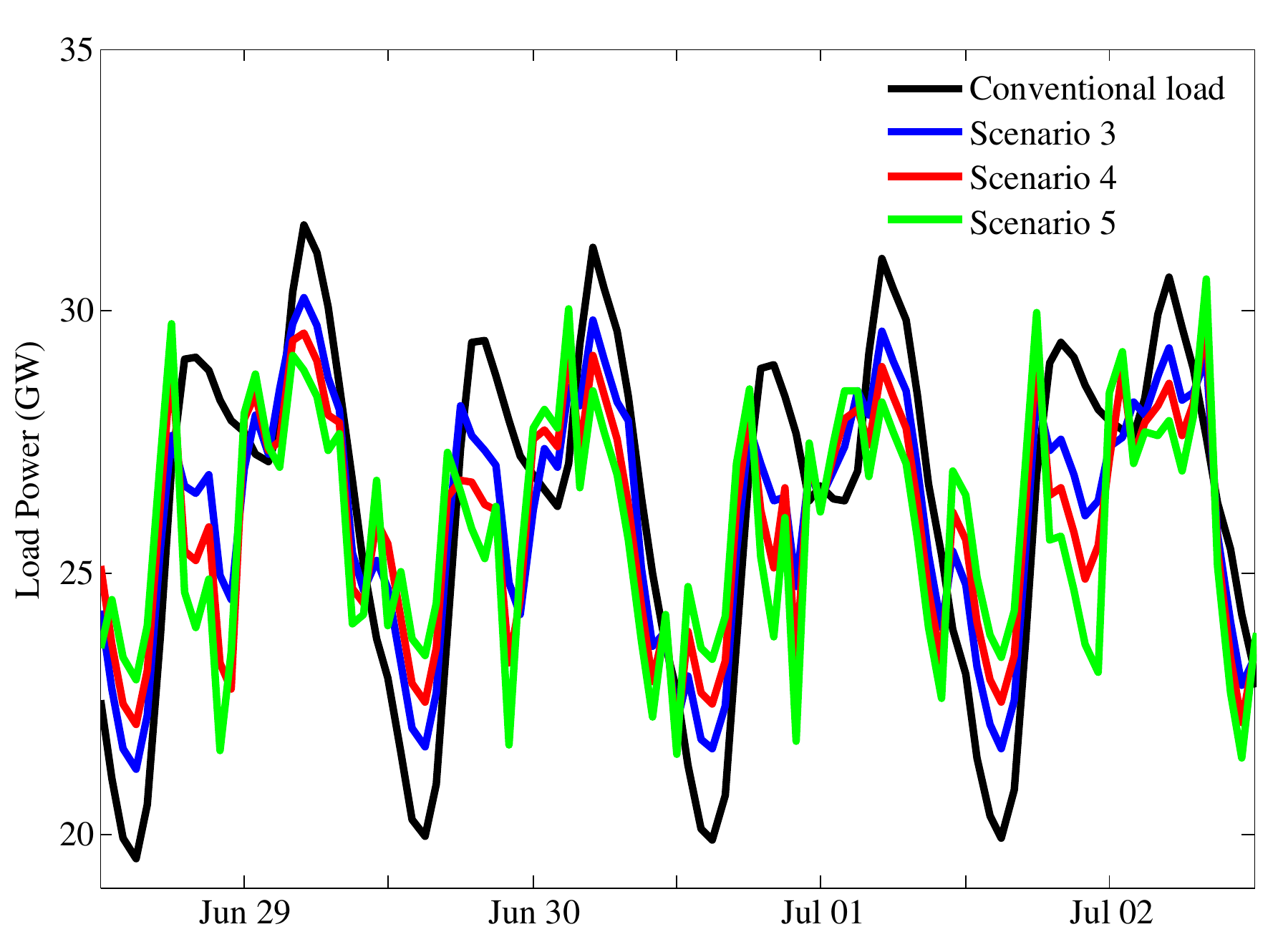}
\label{fig_second_caseA2}}} 
\caption{Demand profile for the NEM for Scenarios 2 to 5 in one of the critical (a) summer and (b) winter peaks}
\label{figure:Demand profile for scenarios 2 to 5 in one of the critical (a) summer and (b) winter peak}
\end{figure}

\begin{figure} 
{\subfloat[]{\includegraphics[width=9.4cm, height=6.5cm]{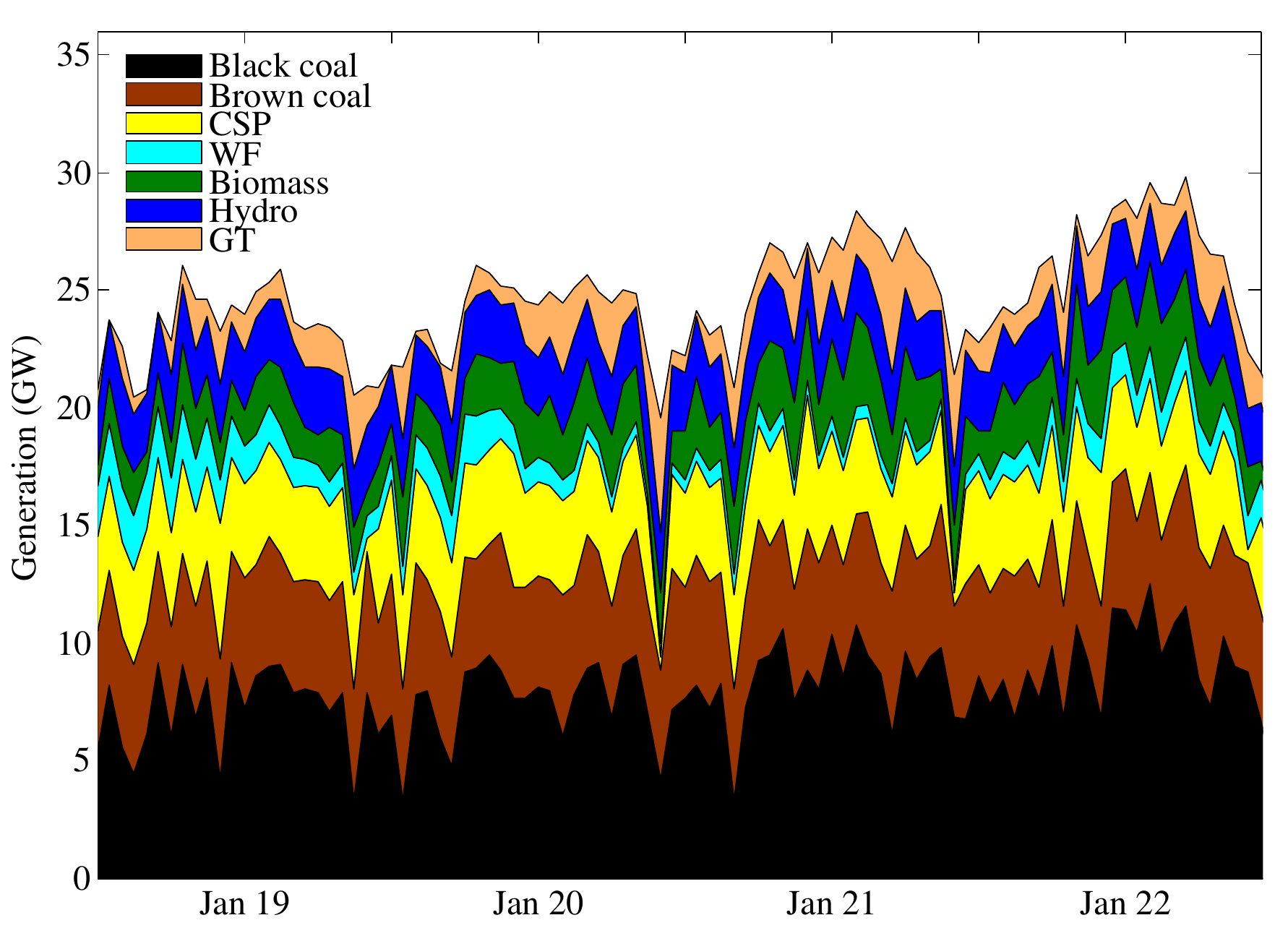}
\label{fig:first_caseA1}}}
{\subfloat[]{\includegraphics[width=9.4cm, height=6.5cm]{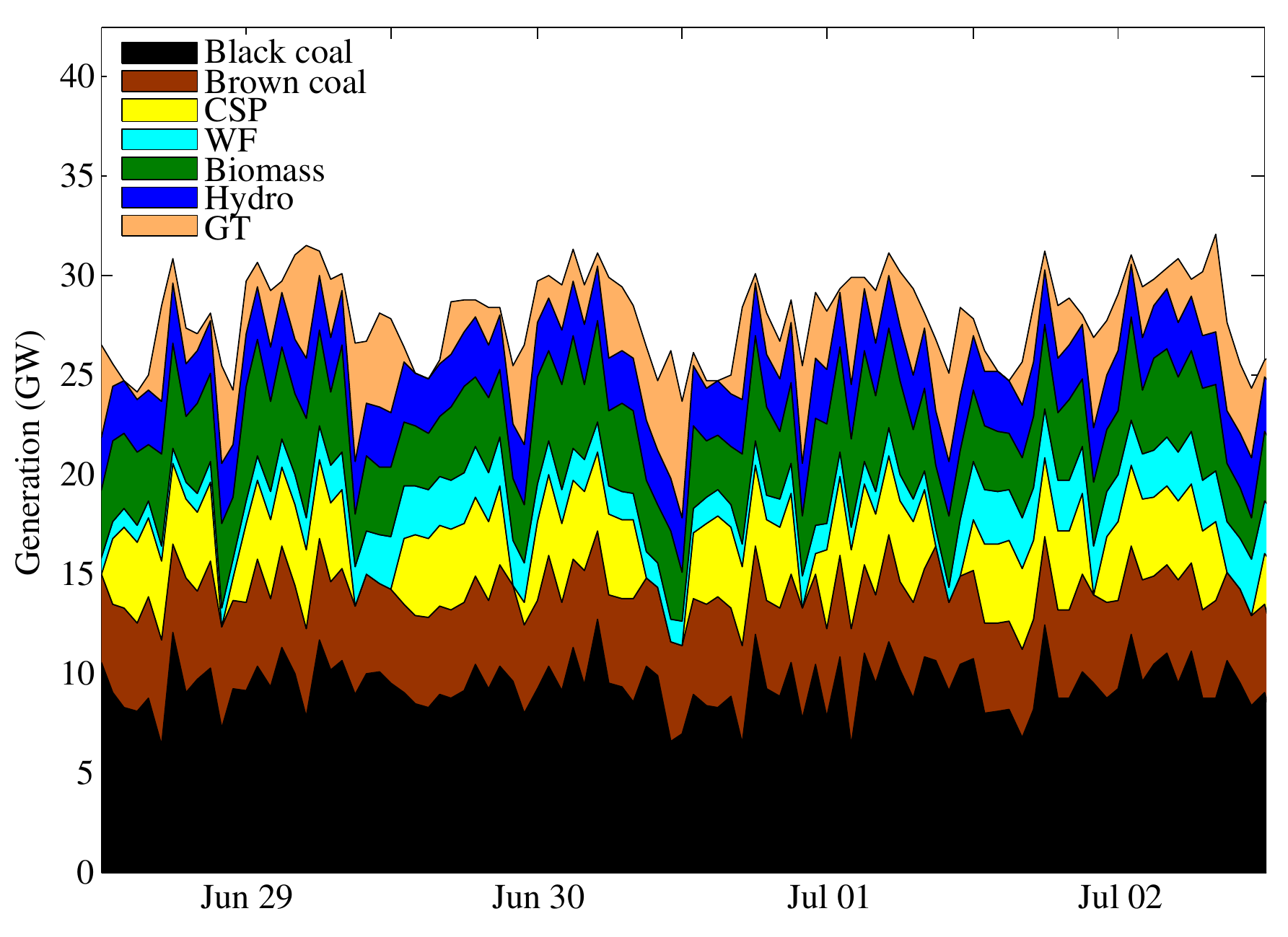}
\label{fig_second_caseA2}}}
\caption{Balancing results for the NEM for the medium DR Scenario in one of the critical (a) summer and (b) winter peaks}
\label{figure:Balancing results for scenario 3 from 10th to 15th of January 2020}
\end{figure}

 Figure~\ref{figure:Demand profile for scenarios 2 to 5 in one of the critical (a) summer and (b) winter peak} shows the conventional load profile and the effect of different uptake scenarios on the load profile in one of the critical summer (i.e. 19\superscript{th} to 22\superscript{nd} of January) and winter peaks (i.e. 29\superscript{th} June to 2\superscript{nd} July) for the NEM. The balancing results for the medium DR scenario during those peak hours are demonstrated in Figure~\ref{figure:Balancing results for scenario 3 from 10th to 15th of January 2020}. As shown in Figure~\ref{figure:Balancing results for scenario 3 from 10th to 15th of January 2020}, in the critical summer days the wind is not strong enough and the output of WF is low, and in critical winter days the solar exposure reduces and CSP output decreases as well. However, due to enough response capacity from price-responsive users (Figure~\ref{figure:Demand profile for scenarios 2 to 5 in one of the critical (a) summer and (b) winter peak}) and backup supply, balancing under these worst-case conditions are maintained. Figures~\ref{figure:Demand profile for scenarios 2 to 5 in one of the critical (a) summer and (b) winter peak} and ~\ref{figure:Balancing results for scenario 3 from 10th to 15th of January 2020} show that management of demand-side technologies can help balance fluctuating RESs power and demand in real-time even under the worst-case conditions. It also improves the balancing and loadability with the increased penetration of RESs in the network.

\section{Conclusion}

\par In this paper, an aggregate load model considering demand-side technologies is proposed for FG scenarios. The load model is intended to be used for system studies at transmission levels. The model is inspired by the smart home concept and formulated as an optimization problem aiming at minimizing the electricity cost. The users are assumed to be price takers. Also, the effect of the load model on the load profile, balancing and loadability of the NEM with the increased penetration of RESs is studied using a modified 14-generator model of the NEM. 

\par Simulation results show that with the increased penetration of RESs and no price-responsive users the loadability is reduced. With price-responsive users, however, the loadability is improved and the required backup supply is reduced with increasing uptake of PV-plus-storage and DR. Interestingly, increasing the penetration of demand-side technologies beyond a certain point does not necessarily improve the performance further and might even deteriorate the system loadability. This depends on the price taking assumption of the loads, when secondary peaks are created due to load synchronization. The proposed model is the first step made towards demand modelling (including DG, storage and DR) for FG scenarios. A more realistic analysis of high penetration of demand-side sources and storage will require an implicit modelling of demand response aggregators using game-theoretic approaches, which will be the focus of future research. Also, the effect of the proposed modelling approach will be studied on power system stability.


\begin{thebibliography}{16}

\bibitem{DepartofResources}
Australian Government, Department of Industry, ``Energy Facts, Statistics and Publications", (Department of Industry of Australian Government, 2014). [Online]. Available: http://www.innovation.gov.au/Energy/Pages/default.aspx.

\bibitem{EPRI}
EPRI, ``The Integrated Grid Realizing the Full Value of Central and Distributed Energy Resources", (EPRI, 2014), pp. 1--44.

\bibitem{RockyMountain}
Rocky Mountain Institute,  Homer Energy, Cohnreznick Think Energy, ``The Economics of Grid Defection When and Where Distributed Solar Generation Plus Storage Competes with Traditional Utility Service", Tech. Rep., 2014, pp. 1--73.

\bibitem{PV}
Parkinson,~G.: ``People power: Rooftop solar PV reaches 3GW in Australia", 2013. [Online]. Available: http://reneweconomy.com.au/2013/people-power-rooftop-solar-pv-reaches-3gw-in-australia-99543.

\bibitem{AEMO2012A}
AEMO, ``2012 NTNDP Assumptions and Inputs", 2012. [Online]. Available: http://www.aemo.com.au/Electricity/Planning/National-Transmission-Network-Development-Plan/Assumptions-and-Inputs.

\bibitem{AEMO2012}
AEMO, ``Solar PV Forecast for AEMO 2012-2022", (SunWiz and SolarBusinessServices, 2012), pp. 1--49.

\bibitem{CSIRO}
CSIRO FG Forum, ``Future Grid Forum: change and choice for Australia electricity system", 2013. [Online]. Available:http://www.csiro.au/Organisation-Structure/Flagships/Energy-Flagship/Future-Grid-Forum-brochure.aspx.

\bibitem{ATA}
Szatow,~T., Moyse,~D.: ``What Happens When We Un-Plug? Exploring the Consumer and Market Implications of Viable, off-Grid Energy Supply, Research Phase 1: Identifying off-Grid Tipping Points", (Energy for the People and Alternative Technology Association, 2014), pp. 1--58.

\bibitem{Energy2010}
Wright,~M., Hearps,~P.: ``Zero Carbon Australia Stationary Energy Plan", (The University of Melbourne Energy Research Institute, 2010), pp. 1--171.

\bibitem{Elliston2012}
Elliston,~B., Diesendorf,~M., MacGill,~I.: ``Simulations of Scenarios with 100\% Renewable Electricity in the Australian National Electricity Market", Energy Policy, 2012, 45, pp. 606--613.

\bibitem{Elliston2013}
Elliston,~B., MacGill,~I., Diesendorf,~M.: ``Least Cost 100\% Renewable Electricity Scenarios in the Australian National Electricity Market", Energy Policy, 2013, 59, pp. 270--282.

\bibitem{Budischak2013}
Budischak,~C., Sewell,~D., Thomson,~H., Mach,~L., Veron,~D.~E., Kempton,~W.: ``Cost-Minimized Combinations of Wind Power, Solar Power and Electrochemical Storage, Powering the Grid up to 99.9\% of the Time", Journal of Power Sources, 2013, 225, pp. 60--74.

\bibitem{Hart2011}
Hart,~E.~K., Jacobson,~M.~Z.: ``A Monte Carlo Approach to Generator Portfolio Planning and Carbon Emissions Assessments of Systems with Large Penetrations of Variable Renewables", Renewable Energy, 2011, 36, (8), pp. 2278--2286.

\bibitem{Tischer2011}
Tischer,~H., Verbic,~G.: ``Towards a Smart Home Energy Management System-A Dynamic Programming Approach", in the Innovative Smart Grid Technologies Asia, 2011.
\bibitem{Tsang2006}
Tsang,~I.~W.-H., Kwok,~J.~T.-Y., Zurada,~J.~M.: ``Generalized Core Vector Machines," IEEE transactions on neural networks / a publication of the IEEE Neural Networks Council, 2006, 17, (5), pp. 1126--40.

\bibitem{Gibbard2010}
Gibbard~M., Vowles,~D.: ``Simplified 14-Generator Model of the SE Australian Power System", (The University of Adelaide, 2010), pp. 1--45.

\bibitem{Kundur}
Kundur, P.: ``Power System Stability and Control" (EPRI Power System Engineering Series, McGraw-Hill, 1994, 1st edn.).

\bibitem{Van}
Van Cutsem, T., Vournas, C.: ``Voltage Stability of Electric Power Systems" (Kluwer International Series in Engineering and Computer Science, 1998, 1st edn.), pp. 93--132.

\bibitem{Samadi2012}
Samadi,~P., Mohsenian-rad,~A.-H., Schober,~R., Wong,~V.W.S.: ``Advanced Demand Side Management for the Future Smart Grid Using Mechanism Design", IEEE Transactions on Smart Grid, 2012, 3, (3), pp. 1170--1180.

\bibitem{Szkuta1999}
Szkuta,~B.R., Sanabria,~L.a., Dillon,~T.S.: ``Electricity Price Short-Term Forecasting using Artificial Neural Networks", IEEE Transactions on Power Systems, 1999, 14, (3), pp. 851--857.


\bibitem{Yamin2004}
Yamin,~H., Shahidehpour,~S., Li,~Z.: ``Adaptive Short-Term Electricity Price Forecasting using Artificial Neural Networks in the Restructured Power Markets", International Journal of Electrical Power \& Energy Systems, 2004, 26, (8), pp. 571--581.


\bibitem{Niu2010}
Niu,~D., Liu,~D,, Wu,~D.~D.: ``A Soft Computing System for Day-Ahead Electricity Price Forecasting", Applied Soft Computing, 2010, 10, (3), pp. 868--875.

\bibitem{Energy2012}
Australian Government, Bureau of Resources and Energy Economics, ``Australian Energy Technology Assessment", Tech. Rep., 2012.

\bibitem{Chapman2013}
Chapman,~A.C., Verbic,~G., Hill,~D.J.: ``A Healthy Dose of Reality for Game-Theoretic Approaches to Residential Demand Response", in Bulk Power System Dynamics and Control - IX Optimization, Security and Control of the Emerging Power Grid (IREP) 2013 IREP Symposium.






\end{thebibliography}
\end{document}